\documentclass[conference]{IEEEtran}
\IEEEoverridecommandlockouts
\usepackage{amsmath,amssymb,amsfonts}
\usepackage{algorithmic}
\usepackage{amsmath}
\usepackage{mathtools}

\usepackage{graphicx}
\usepackage{kotex}
\usepackage{afterpage}
\usepackage{textcomp}
\usepackage{xcolor}
\def\BibTeX{{\rm B\kern-.05em{\sc i\kern-.025em b}\kern-.08em
    T\kern-.1667em\lower.7ex\hbox{E}\kern-.125emX}}
\usepackage[hidelinks]{hyperref} 
\usepackage[
backend=biber,
style=numeric-comp,
   sorting=none
]{biblatex}

\begin{document}

\title{Enhancing ASR Performance through OCR Word Frequency Analysis: Theoretical Foundations}
\author{
Kyudan Jung$^{\star}$, Nam-Joon Kim$^{\dag}$, Hyun Gon Ryu$^{\ddag}$, Hyuk-jae Lee$^{\dag}$\\
$^{\star}$ Dept. of Mathematics, Chung-Ang University, Seoul, Korea \\
$^{\dag}$ Dept. of Electrical and Computer Engineering, Seoul National Univeristy, Seoul, Korea \\
$^{\ddag}$ NVIDIA, Seoul, Korea
}

\maketitle
\begin{abstract}

As the interest in large language models grows, the importance of accuracy in automatic speech recognition has become more pronounced. This is especially true for lectures that include specialized terminology. In such cases, the success rate of traditional ASR models tends to be low, presenting a significant challenge. A method using the word frequency difference approach has been proposed to improve ASR performance for specialized terminology. We investigated this proposal through experiments and data analysis to determine if it effectively addresses the issue. In addition, we introduced the power law as the theoretical foundation for the relative frequency methodology mentioned in this approach. \end{abstract}


\section{Introduction}

Recent advancements in conversational AI and large language model (LLM) research have been notably vigorous.
These advancements extend their applications to everyday conversations and specialized fields like medicine, law, engineering, and education. To facilitate this, Automatic Speech Recognition (ASR) systems must be capable of accurately recognizing and understanding specialized terminology within these fields.
Consequently, there is ongoing research to enhance the recognition rates of specialized terminology. Ma's study implemented a method in which lip movement analysis provides clues to enable ASR systems to recognize unfamiliar words \cite{ref1}. This represents an innovative attempt to utilize visual information to compensate for unclear pronunciation in audio signals. Furthermore, Guo's research involved an approach to correct inaccuracies in ASR outputs by incorporating spelling-correction efforts \cite{ref2}. Guo proposed a method that assigns weights according to spelling rules to the n-best sequences obtained from the ASR model and selects the sequence with the highest score.

Jung's research proposed a method to improve the ASR performance of lecture audios by extracting text from lecture videos using Optical Character Recognition (OCR) \cite{ref3}. The extracted text is then utilized to enhance the ASR performance.

In this study, we experimented with Jung's method, dissecting and analyzing it in detail. The results identified two problems that could be effectively resolved through the development of this method. Additionally, the theoretical background of the Relative Frequency (RF) methodology was validated through the power law.

\begin{figure*}[t]    \centering    \includegraphics[width=\textwidth, height=5cm]{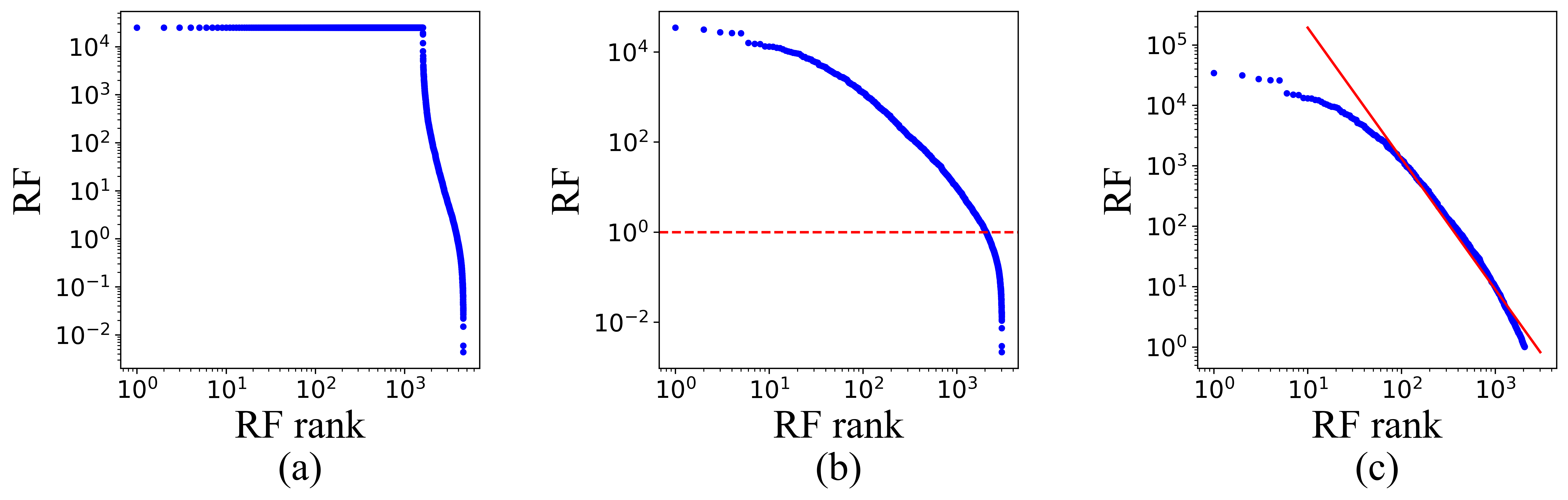}    \caption{(a) Graph of RF using the existing method \cite{ref3}, with log(RF rank) on the x-axis and log(RF) on the y-axis. (b) Graph of RF for data processed using Method 1 in addition to the existing method. (c) Graph of RF for data processed using both Method 1 and Method 2 in addition to the existing method.}  
\end{figure*}

\section{Word Frequency Difference Approach}

Jung's proposed word frequency difference approach is defined as follows. First, the Normal Frequency (NF) of a word, the Frequency of the word in a Lecture (LF), and the Relative Frequency (RF) are defined. All three metrics were determined for a single word. 
\begin{itemize}   
\item 
\textbf{NF}; Normal Frequency: Represents the frequency of a word in general contexts. It is defined in Equation 1 \cite{ref4} as the count of a word divided by the total counts of all words using the Google Web Trillion Word Frequency Dataset \cite{ref5}, which includes approximately 333,000 of the most commonly used words on the Internet and their count values. $c_{\text{NF}}$ in Equation 1 is defined as in Equation 2. A similarly large corpus can also be used instead of the Google Web Trillion dataset. In this paper, we refer to such a large corpus as a large-text dataset (LTD). 
\begin{equation}    \text{NF}_{w_i} \coloneqq \frac{c_{\text{NF}}(w_i)}{\sum_{w_i \in \text{LTD}} c_{\text{NF}}(w_i)}\end{equation}\begin{equation}c_{\text{NF}}(w_i) := \begin{cases} c(w_i) & \text{if } w_i \in \text{LTD} \\\min(\text{count}) & \text{if } w_i \notin \text{LTD}\end{cases}
\end{equation}
where $ w_i $ is the word that make up the output sentence of the beam search decoder, and  $ c(w_i) $ is 'count' value of word in LTD.
\item \textbf{LF}; Lecture Frequency: Represents the frequency of a word in a lecture context as in Equation 3.\cite{ref4} In this methodology, it is defined as the count of the targeted word among all words extracted via OCR, divided by the total number of all words.    
\begin{equation}    LF := \frac{\text{\# of specific words}}{\text{\# of words in RF word pool}}
\end{equation}

\item \textbf{RF}; Relative Frequency: It is defined as in Equation 4.\cite{ref4} Thus, if a word appears more frequently in lectures than in general contexts, the RF value increases. Jung determined that words used more frequently in lectures than in general are likely specialized terms for that context, and enhanced the ASR performance for such terms by making them appear more prominently in the ASR decoder.

\begin{equation}    RF := \frac{\text{LF}}{\text{NF}}\end{equation}
\end{itemize}

In subsequent research, Jung modified the NF data processing module, resulting in an improvement of up to 3.22\%. This enhanced method applied two main approaches:
\begin{itemize}    
\item  
\textbf{Method 1.} When calculating the count value for NF, if a word in the output sentence of the beam search decoder, which is the decoder for ASR, did not exist in the LTD, its count value was set to zero. In the revised method, if a word is absent in the LTD, its count value is replaced by the minimum count value among all the words identified in the OCR. Whereas the original method used the count values of all words in the LTD to compute NF, the revised method calculates NF only for the set of words extracted using OCR, referred to as the OCR dataset (OCRd).    
\item \textbf{Method 2.} To create a modified RF that follows the Power Law, all RF values less than 1 are replaced with 1 in the data, as in Equation 5. 

\begin{equation}
    RF := 
\begin{cases} 
RF_{\text{max}} & \text{if } RF \geq 1 \\
1 & \text{otherwise}
\end{cases}
\end{equation}

\end{itemize}

After calculating RF in this manner, the OCR score is computed and then combined with the score from the existing beam search decoder. Based on this Equation 6,
\begin{equation}
    RF(w) \propto \text{rank}_{RF}(w)^{k}
\end{equation}
where w is the word, and $\text{rank}_{RF}$ is the rank of the RF of 
the words in descending order, $\text{score}_{OCR}$ is defined by Equation 7.

\begin{equation}
    \text{score}_\text{OCR} = \frac{1}{m} \sum_{i=1}^{m}{1-RF(w_i)^{\frac{1}{k}}}
\end{equation}
where m is the total number of words in the output 
sentence of the beam search decode, $RF(w_i)$ is the RF value 
of the $i$th word, and k is the same value as in Equation 1, which is 
the slope of the straight line with a linear regression of the 
log-log scale rank of RF-RF graph.

Finally, $\text{score}_{\text{OCR}}$ is merged with $\text{score}_{\text{old}}$ which comes from original ASR models' beam search decoder.
\begin{equation}
    \text{score}_{\text{new}} = \text{score}_{\text{old}} + \lambda_{OCR} \cdot \text{score}_{\text{OCR}}
\end{equation}
where lambda is a hyperparameter with default value of 1.0. 

\section{Analysis and Experiments}
This section introduces the methodologies and experiments designed to address discrepancies between text extracted by OCR and a large-scale language corpus. We analyze the issues arising when OCR-extracted words are not found in the large text database (LTD), comparing the performance of existing and improved methods.

\subsection{Method 1. Handling Words Extracted by OCR Not Found in the LTD }
When a word extracted by OCR is not found in LTD, the existing method treats its count value as zero. This might seem justified because the word is not present in the LTD, but it results in an NF of zero, making it impossible to define the RF.

Additionally, the minimum word count value in the existing Google Corpus is 12,711, which creates a significant discrepancy with the other word counts. To address this, \cite{ref3} suggested replacing the RF values of such words with the maximum RF value of the current frame. We conducted the following experiment to understand the drawbacks of existing methods. 

\subsubsection{Experiment on the Existing Method}We used 108 hours of lecture videos from the course "Advanced Compiler" with the LTD and OCRd datasets. LTD is the Google Web Trillion Word Frequency Dataset, and OCRd is a dataset of all words captured at 30-second intervals from lecture videos, predominantly consisting of professors’ PowerPoint slides. These two datasets were merged based on common words in the 'word' column. This process resulted in a dataset of 334,935 words (rows), combining 333,333 words from LTD and 1,601 words recognized by OCR, but not LTD.

For OCRd words that were not found in LTD, their counts were set to zero. NF was calculated as originally defined by dividing each count by the total sum of counts. When NF was zero, RF was replaced with the maximum RF value, reasoning that if a word's NF is zero, it is likely to be a specialized term not present in the LTD.

\subsubsection{Data Analysis of the Experiment on the Existing Method}Fig.1a shows the graph based on the 334,935 words mentioned above, with log(RF rank) on the x-axis and log(RF) on the y-axis. The existing method of calculating RF shows that RF values maintain a constant range, equalizing the RF scores of top-ranked words and reducing the influence of terms specifically used in lectures compared to general contexts.

This confirmed that the existing method reduces RF reliability and accuracy by uniformly assigning high RF values to words not found in the LTD.

\subsubsection{Experiment on the Improved Method 1}
When applying the improved data processing method proposed in \cite{ref4}, enhanced RF values were observed. The datasets from both LTD and OCRd were merged based on the intersection of words in the 'word' column of the latter. This process excludes words from the LTD that are not extracted by the OCR. When a word extracted by OCR is not in LTD, its count value is replaced by the minimum count value from the merged dataset. The resulting graph appears in Fig.1b.

\subsection{Method 2: Observing Compliance with Power Law for Words with RF Values Above 1}Fig.1b shows the RF values for all words calculated using the method proposed in \cite{ref4} for LTD and OCRd texts. The graph is on a log-log scale with the RF rank on the x-axis and RF values on the y-axis. According to Zipf's law, the word frequency in a text collection roughly follows a power law. However, words with RF values less than one showed a different trend.

In Fig.1b, the red horizontal dotted line represents 1. The linear approximation had 4545 degrees of freedom and a residual standard error (RSE) of 0.6606.Fig.1c shows the graph drawn excluding data with an RF less than 1. These data had a linear approximation with 3583 degrees of freedom and an RSE of 0.3401, leading to the observation that the rank of the RF-RF data follows the power law well.

\section{Conclusion}In this paper, we demonstrated the validity of the proposed method as an improvement over the previous approach. Specifically, we showed through experimentation that the approach of has drawbacks in three aspects and that the method in can resolve these shortcomings. Notably, the improved method has a convincing theoretical foundation based on the power law.

\printbibliography

\end{document}